\documentclass[12pt]{amsart}

\usepackage{amsmath}
\usepackage{amssymb}
\usepackage{multirow}
\usepackage{epsfig}

\usepackage{pst-node}
\usepackage{tikz}
\usepackage{hyperref}

\usepackage{a4wide}

\renewcommand{\proof}{\par\noindent{\it Proof.\ \ }}
\def\qed{\ifmmode\square\else\nolinebreak\hfill
$\square$\fi\par\vskip12pt}

\newtheorem{propo}{Proposition}[section]
\newtheorem{defi}[propo]{Definition}

\newtheorem{lemma}[propo]{Lemma}

\newtheorem{theo}[propo]{Theorem}
\newtheorem{examp}[propo]{Example}

\newtheorem{prob}[propo]{Problem}

\newtheorem{rem}[propo]{Remark}
\newtheorem{constr}[propo]{Construction}

\newcommand{\bl}{\begin{lemma}}
\newcommand{\el}{\end{lemma}}

\def\HA{{\rm HA}}
\def\HS{{\rm HS}}

\def\AS{{\rm AS}}
\def\TW{{\rm TW}}
\def\SD{{\rm SD}}

\def\PA{{\rm PA}}

\def\Out{{\rm Out}}

\def\K{{\rm K}}

\usepackage{indentfirst,latexsym,bm}

\def\Sym{{\rm Sym}}

\def\Aut{{\rm Aut}}
\def\Out{{\rm Out}}

\def\diam{{\rm diam}}

\def\H{{\rm H}}
\def\Cos{{\rm Cos}}

\begin{document}
\title{ Vertex quasiprimitive two-geodesic transitive graphs }

\thanks{Supported by the NNSF of China (12061034,12071484) and  NSF of Jiangxi (20192ACBL21007, GJJ190273)}

\thanks{The author is grateful to Professors Cheryl Praeger, Cai Heng Li and Alice Devillers for their discussion
and comments on this paper.}

\author[W. Jin]{Wei Jin}
 \address{Wei Jin\\School of Mathematics and Statistics\\
Central South University\\
Changsha, Hunan, 410075, P.R.China}
\address{School of Statistics\\
Jiangxi University of Finance and Economics\\
 Nanchang, Jiangxi, 330013, P.R.China}
\email{jinweipei82@163.com}
%\author[W. J. Liu]{Wei Jun Liu}
%\address{Wei Jun Liu\\School of Mathematics and Statistics\\
%Central South University\\
%Changsha, Hunan, 410075, P.R.China}
%\email{wjliu@csu.edu.cn}

%\date\today

\maketitle

\begin{abstract}

For a non-complete  graph $\Gamma$,   a vertex triple $(u,v,w)$ with
$v$ adjacent to both $u$ and $w$ is called a $2$-geodesic if
$u\neq w$ and $u,w$ are not adjacent.  Then $\Gamma$ is said to be
$2$-geodesic transitive if its automorphism group is
transitive on both arcs and 2-geodesics. In previous work  the author showed that if a $2$-geodesic transitive graph $\Gamma$ is locally disconnected and its automorphism group $\Aut(\Gamma)$ has a non-trivial normal subgroup
which is intransitive on the vertex set of $\Gamma$, then $\Gamma$ is a cover of a smaller 2-geodesic transitive graph. Thus the `basic' graphs to study are
those for which $\Aut(\Gamma)$ acts quasiprimitively on the vertex set.
In this paper, we study  2-geodesic transitive graphs which are locally disconnected and $\Aut(\Gamma)$ acts quasiprimitively on the vertex set.
We first determine all the possible quasiprimitive action types and give examples for them,
and then classify the family of  $2$-geodesic transitive graphs  whose automorphism group  is primitive on its vertex set  of $\PA$ type.

%This paper partly answers \cite[Question 8.4]{GLP1} for $s=3$.

\end{abstract}

\vspace{2mm}

%\hspace{-17pt}{\bf AMS:} 62P05

%\hspace{-17pt}{\bf Subject category:}IM10

 \hspace{-17pt}{\bf Keywords:}  $2$-geodesic transitive graph, quasiprimitive permutation group, automorphism group.

 \hspace{-17pt}{\bf Math. Subj. Class.:} 05E18; 20B25

\section{Introduction}

A vertex triple $(u,v,w)$ in a non-complete
graph $\Gamma$ with $v$ adjacent to both $u$ and $w$ is a
\emph{$2$-arc} if $u\neq w$, and a \emph{$2$-geodesic} if $u,w$ are
not adjacent. A graph $\Gamma$ is said to be  \emph{$(G,2)$-geodesic
transitive} or \emph{$(G,2)$-arc transitive} if, the group  $G$ is
transitive on  arc set, and also transitive on 2-geodesic set or 2-arc set, respectively,  where  $G\leq \Aut(\Gamma)$.
If $G= \Aut(\Gamma)$,  then $\Gamma$ is simply said to be  \emph{$2$-geodesic transitive} or \emph{$2$-arc transitive}.   Clearly, every 2-geodesic is a 2-arc, but some 2-arcs may not be 2-geodesics. If $\Gamma$ has girth 3 (length of the
shortest cycle is 3), then the 2-arcs contained in $3$-cycles are
not 2-geodesics. The graph in Figure 1 is the Kneser graph $KG_{6,2}$ which is
$2$-geodesic transitive but not $2$-arc transitive with valency 6.
Thus the family of non-complete $2$-arc transitive graphs is
properly contained in the family of $2$-geodesic transitive graphs.   Devillers, Li,  Praeger and the author in \cite[Theorem 1]{DJLP-clique} proved that if
$\Gamma$ is  $(G,2)$-geodesic transitive of valency at least 2, then for each vertex
$u$,  either

{\rm (1)}  $[\Gamma(u)]$ is connected of diameter $2$; or

{\rm (2)}  $[\Gamma(u)]\cong
m\K_r$ for some  integers $m\geq 2,r\geq 1$.

\begin{figure}[t]
\centering
\includegraphics[width=4.2cm,height=4.2cm]{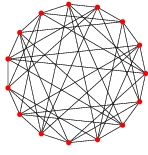}
\caption{Kneser graph $KG_{6,2}$} \label{kneser}
\end{figure}

The graphs in case (1) were studied
in \cite{J-con-2017}, and such  graphs which are locally isomorphic to $\overline{mC_n}$ were completely determined.
In \cite{JT-local}, the author investigated the case (2), and they  gave a reduction theorem for such graphs.
In this paper, we
continue the study  of locally disconnected $(G,2)$-geodesic transitive graphs.

We denote by $\mathcal{F}(m,r)$ the family of connected graphs $\Gamma$ such that
$[\Gamma(u)]\cong m\K_r$ for each vertex $u$, and for fixed $m\geq
2,r\geq 1$.
Let  $\Gamma \in \mathcal{F}(m,r)$ be a $(G,2)$-geodesic transitive
graph where  $m,r\geq 2$. If $G$ is not quasiprimitive on
the vertex set, then \cite{JT-local} showed that  there exists an intransitive normal subgroup $N$
of $G$ such that: $\Gamma$ is a cover of $\Gamma_N$, $G/N$ is
quasiprimitive on $V(\Gamma_N)$, and  either $\Gamma_N\cong
\K_{mr+1}$ is $G/N$-arc transitive or  $\Gamma_N\in \mathcal{F}(m,r)$
is non-complete $(G/N,2)$-geodesic transitive.

We study `basic' locally disconnected $(G,2)$-geodesic transitive graphs,
that is,   we suppose that $G$ is quasiprimitive on the vertex set.
Our main theorem determines all the possible quasiprimitive action types.

\begin{theo}\label{2gt-cliq-basic-2}
Let $\Gamma \in \mathcal{F}(m,r)$  be a $(G,2)$-geodesic transitive graph  where $m,r\geq 2$ and $G\leq \Aut(\Gamma)$. Suppose that $G$ is
quasiprimitive on $V(\Gamma)$ of type $X$.  Then $X$ is in the set
$\{\HA,\PA,\AS,$ $\TW,\SD\}$.

Further, if the induced action of $G$ on $V(C(\Gamma))$ is not
quasiprimitive, then $X$ is one of
$\{\HA,\PA,\AS,\TW\}$, see Table \ref{Star}; if the induced action of $G$ on $V(C(\Gamma))$ is  quasiprimitive of type
$Y$, then $(X,Y)$ is one of
$\{(\HA,\HA),(\PA,\PA),(\AS,\AS),$ $(\TW,\TW),(\SD,\PA),(\PA,\SD)\}$, see Table \ref{Biprimitive}.
\end{theo}

\begin{table}[t]\caption{Quasiprimitive of type $X$ on  vertices  but not on maximal cliques }\label{Star}
\medskip
\centering
\begin{tabular}{|c|c|c|c|c|c|}
\hline Case     &   & Example    \\
\hline   $X$    & $\HA$ &  Examples \ref{2gt-star-ha} and  \ref{2gt-qp-examphamm12}     \\
\hline  $X$    & $\PA$ &   Example \ref {2gt-qp-examphamm1} \\
\hline $X$    & $\AS$   &   Examples \ref{subdivision-as-star-1} and  \ref{2gt-star-as}  \\
\hline $X$    & $\TW$   &   We don't have examples. \\
  \hline
\end{tabular}
\end{table}

\begin{table}[]\caption{Quasiprimitive of types $X$ on  vertices and $Y$ on maximal cliques }\label{Biprimitive}
\medskip
\centering
\begin{tabular}{|c|c|c|c|c|c|}
\hline Case     &  $Y$ & Example    \\
\hline   $X$    & $(\HA,\HA)$ &   We don't have examples.   \\
\hline  $X$    & $(\PA,\PA)$ &    We don't have examples. \\
\hline  $X$    & $(\AS,\AS)$   & Examples \ref{smallval6-psl-1}, \ref{2gt-qp-examp1} and \ref{primitive-lemma-1}   \\
\hline $X$    & $(\SD,\PA)$   &   Example \ref{2gt-cliq-sd-pa} \\
\hline $X$    & $(\PA,\SD)$   &  Example \ref{2gt-cliq-sd-pa} \\
\hline $X$    & $(\TW,\TW)$   &  We don't have examples.  \\
  \hline
\end{tabular}
\end{table}

We give a remark of Theorem \ref{2gt-cliq-basic-2}.

\begin{rem}\label{pls-r1-1}
{\rm  (1) A  $(G,2)$-geodesic transitive graph $\Gamma \in
\mathcal{F}(m,1)$  has girth at least $4$, and hence each $2$-arc is
a $2$-geodesic. It follows that $\Gamma$ is  $(G,2)$-arc transitive.
Such graphs have been studied extensively, see
\cite{Baddeley-1,DDP-1,IP-1,Li-Seress-1,Tutte-1,Tutte-2,weiss}.

(2) Let $\Gamma\in \mathcal{F}(m,1)$ for some $m\geq 2$ be
$(G,2)$-geodesic transitive.  Suppose that  $G$ is quasiprimitive on $V(\Gamma)$. Then by Lemmas 5.2 and
5.3 of \cite{Praeger-1993-onanscott}, there are only four
quasiprimitive types  $\{\HA,\TW,\AS,\PA\}$, and examples exist for each
of them.

(3) Let $\Gamma\in \mathcal{F}(m,r)$, $m\geq 2,r\geq 1$. Then by \cite[Theorem 1.4]{DJLP-clique}, $C(\Gamma)\in \mathcal{F}(r+1,m-1)$.
Suppose that $G(\leq \Aut(\Gamma))$ is quasiprimitive on both $V(\Gamma)$ and
$V(C(\Gamma))$ with a minimal normal subgroup  $N$. Suppose that $N$ has a subgroup which is regular on  both
$V(\Gamma)$ and
$V(C(\Gamma))$. Then $|V(\Gamma)|=|V(C(\Gamma))|$, and also
$m=r+1$.  Hence, if $G$ is quasiprimitive on both $V(\Gamma)$ and
$V(C(\Gamma))$ of type $(\HA,\HA)$ or $ (\TW,\TW)$, then $m=r+1$. If $r=1$, then $m=2$, and $\Gamma$ is a cycle.

(4) Let $\Gamma \in \mathcal{F}(m,r)$  be a $(G,2)$-geodesic transitive graph  where $m,r\geq 2$. Suppose that $G$ is quasiprimitive
on $V(\Gamma)$ of type X but not quasiprimitive  on $V(C(\Gamma))$. If X=PA,  then $\Gamma$ is a graph of  Lemma \ref{2gt-star-pa}; if X=\HA, then
$\Gamma$ is the point graph of a partial linear space of  Example \ref{2gt-star-ha}; if X=AS, then
the socle of $G$ and the number of orbits of $G$ on $V(C(\Gamma))$ are in Example \ref{2gt-star-as}.

}
\end{rem}

Examples of
locally disconnected $(G,2)$-geodesic transitive graphs where $G$ is quasi-primitive of type HA, PA and AS will be provided in Examples  \ref{subdivision-as-star-1},  \ref{2gt-star-as}, \ref {2gt-qp-examphamm1}, and  Examples \ref{2gt-star-ha}, \ref{2gt-qp-examphamm12}. While  Examples  \ref{smallval6-psl-1}, \ref{2gt-qp-examp1} and \ref{primitive-lemma-1}
give examples where $G$ is quasiprimitive of type $\AS$ on  both vertices and  maximal cliques.
Example \ref{2gt-cliq-sd-pa} gives examples where $G$ is quasiprimitive of type $\SD$ on   vertices and $\PA$ on maximal cliques and also gives
examples where $G$ is quasiprimitive of type $\PA$ on   vertices and $\SD$ on maximal cliques.
We do not investigate the case where $G$ is quasiprimitive of type TW on vertices and  the cases where $G$ is quasiprimitive of type HA or PA  on both vertices and  maximal cliques.

Our second theorem determines precisely
the class of  $2$-geodesic transitive graphs  whose automorphism group  is primitive on its vertex set  of $\PA$ type.

\begin{theo}\label{2gt-basic-th2}
Let $\Gamma$ be a $2$-geodesic transitive graph. Then  $\Aut(\Gamma)$ is primitive on $V(\Gamma)$ of $\PA$ type if and only if   $\Gamma$ is either the Hamming graph $\H(d,n)$ for some  $d\geq 2$ and $n\geq 5$ or $\overline{\H(2,n)}$ for some  $n\geq 5$.
\end{theo}

Let $\Gamma \in \mathcal{F}(m,r)$  be a $(G,2)$-geodesic transitive graph  where $m,r\geq 2$. Suppose that $G$ is
quasiprimitive on $V(\Gamma)$.
To finish the classification of  such graphs, we pose the following two problems.

\begin{prob}\label{}
Let $\Gamma \in \mathcal{F}(m,r)$  be a $(G,2)$-geodesic transitive graph  where $m,r\geq 2$. Suppose that $G$ is
quasiprimitive on $V(\Gamma)$ of type $X$.

{\rm (1)} If the induced action of $G$ on $V(C(\Gamma))$ is not
quasiprimitive, then give examples for $X=\TW$.

{\rm (2)} If the induced action of $G$ on $V(C(\Gamma))$ is  quasiprimitive of type
$Y$, then give examples for $(X,Y)=(\HA,\HA),(\PA,\PA)$ or $(\TW,\TW)$.

\end{prob}

\begin{prob}\label{}
Let $\Gamma \in \mathcal{F}(m,r)$  be a $(G,2)$-geodesic transitive graph  where $m,r\geq 2$. Suppose that $G$ is
quasiprimitive on $V(\Gamma)$ of type $X$.

{\rm (1)} If X=\TW and  the induced action of $G$ on $V(C(\Gamma))$ is not
quasiprimitive, then  classify such graphs.

{\rm (2)} If the induced action of $G$ on $V(C(\Gamma))$ is  quasiprimitive of type
$Y$, then for each $(X,Y)\in \{(\HA,\HA),(\PA,\PA),(\AS,\AS),$ $(\TW,\TW),(\SD,\PA),(\PA,\SD)\}$, classify such graphs.

\end{prob}

\section{Preliminaries}

In this section, we give some definitions concerning groups, graphs and geometries and also some  results which will be used in our
analysis.

In this paper, all graphs  are finite, simple, connected and undirected. For a graph $\Gamma$, we use  $V(\Gamma)$ and
$\Aut(\Gamma)$ to denote its \emph{vertex set}  and
\emph{ automorphism group}, respectively. For the group theoretic terminology not defined here we refer the reader to \cite{Cameron-1,DM-1,Wielandt-book}.

A transitive permutation group $G$ is said to be \emph{quasiprimitive}, if every non-trivial normal subgroup of $G$ is transitive.
This is a generalization
of primitivity as every normal subgroup of a primitive group is transitive,
but there exist quasiprimitive groups which are not primitive. For knowledge of quasiprimitive permutation groups, see \cite{Praeger-1993-onanscott} and \cite{Praeger-2}.
Praeger \cite{Praeger-1993-onanscott} generalized the O'Nan-Scott Theorem for primitive groups to quasiprimitive
groups and showed that a finite quasiprimitive group is one of eight distinct types: Holomorph Affine (HA), Almost Simple (AS), Twisted Wreath product (TW), Product Action (PA), Simple Diagonal (SD), Holomorph Simple (HS), Holomorph Compound (HC) and Compound Diagonal (CD).

For a subset  $U$ of the vertex set of $\Gamma$, we denote by $[U]$  the
subgraph of $\Gamma$ induced by $U$, and $[\Gamma(u)]$ is the subgraph induced by the neighborhood of the vertex $u$.

 For a graph $\Gamma$, its
\emph{complement} $\overline{\Gamma}$ is the graph with vertex set $V(\Gamma)$,
and two vertices are adjacent if and only if they are not adjacent in $\Gamma$.

An \emph{arc} of  a graph  is an ordered vertex pair such that the two
vertices are adjacent.

For  vertices $v_1$ and $v_2$ in
$V(\Gamma)$,  the smallest value for $n$
such that there is a path of length $n$ from $v_1$ to $v_2$ is
called the \emph{distance} from $v_1$ to $v_2$ and is denoted by
$d_{\Gamma}(v_1, v_2)$. The \emph{diameter} $\diam(\Gamma)$ of
a connected graph $\Gamma$ is the maximum of $d_{\Gamma}(v_1, v_2)$
over all $v_1, v_2 \in V(\Gamma)$. Let $1\leq i\leq \diam(\Gamma)$. For each vertex $v\in V(\Gamma)$, set
$\Gamma_i(v)=\{u|d_{\Gamma}(v, u)=i,u\in V(\Gamma) \}$.

For two vertices $u,v \in
V(\Gamma)$, a \emph{geodesic} from $u$ to $v$ is one of the shortest
paths from $u$ to $v$, and
this  geodesic is an \emph{$i$-geodesic} if $d_{\Gamma}(u,v)=i$.
Suppose that  $G\leq \Aut(\Gamma)$ and  $1\leq s \leq
\diam(\Gamma)$. Then $\Gamma$ is said to be \emph{$(G,s)$-geodesic transitive}
if, for  each $i=1,2,\ldots,s$, $G$ is transitive on the set of $i$-geodesics. When $s=\diam(\Gamma)$,
$\Gamma$ is said to be \emph{$G$-geodesic transitive}.
If  $G=\Aut(\Gamma)$,  we will say that $\Gamma$ is \emph{$s$-geodesic
transitive} or \emph{geodesic transitive}, respectively.

For a positive integer $s$, an $s$-arc of $\Gamma$ is a sequence of
vertices $(v_0,v_1,\ldots,v_s)$  in  $\Gamma$ such that
$v_i,v_{i+1}$ are adjacent and $v_{j-1}\neq v_{j+1}$ where $0\leq
i\leq s-1$ and $1\leq j\leq s-1$. Let $G\leq \Aut(\Gamma)$.  The graph $\Gamma$ is said to be \emph{locally $(G,s)$-arc transitive}
if, for each $i\leq s$, the
stabilizer group $G_u$  is transitive on the set of  $i$-arcs  starting from $u$.  A \emph{ $(G,s)$-arc transitive} graph is a
locally $(G,s)$-arc transitive graph which is also $G$-vertex transitive.

A \emph{clique} of  $\Gamma$ is a complete subgraph and a
\emph{maximal clique} is a clique which is not contained in a larger
clique.   The \emph{clique graph} $C(\Gamma)$ of $\Gamma$ is the
graph with vertex set $\{$all maximal cliques of $\Gamma$$\}$, and
two maximal cliques  are adjacent in $C(\Gamma)$ if and only if they
have at least one common vertex in $\Gamma$.
Let $\Gamma\in \mathcal{F}(m,r)$, $m\geq 2,r\geq 1$. Then by \cite[Theorem 1.4]{DJLP-clique}, $C(\Gamma)\in \mathcal{F}(r+1,m-1)$ and $\Gamma\cong C(C(\Gamma))$. Suppose that $G\leq \Aut(\Gamma)$. Then $G$ induces a faithful action on  $V(C(\Gamma))$.

An \emph{incidence structure} is a triple
$\mathcal{S}=(\mathcal{P},\mathcal{L},\mathcal{I})$ where
$\mathcal{P}$ is a set of \emph{points}, $\mathcal{L}$ is a set of
\emph{lines} and $\mathcal{I}\subseteq \mathcal{P}\times
\mathcal{L}$. If $(p,l)\in
\mathcal{I}$, then we say that the point $p$ and the line $l$ are
incident.  The \emph{incidence graph} of an incidence structure
$\mathcal{S}=(\mathcal{P},\mathcal{L},\mathcal{I})$ is the bipartite
graph with vertex set $\mathcal{P}\cup \mathcal{L}$ and edge set
$\{(p,l)|p\in \mathcal{P},l\in \mathcal{L},(p,l)\in \mathcal{I}\}$.
A \emph{partial linear space} is an incidence structure
$\mathcal{S}=(\mathcal{P},\mathcal{L},\mathcal{I})$ such that any
line is at least incident with two points and any pair of distinct
points is incident with at most one line. In particular,  we say that a partial linear
space has \emph{order} $(m,n)$ if    each point is incident with $m$
lines, and each line is incident with $n$ points.

For an incidence structure $\mathcal{S}=(\mathcal{P},\mathcal{L},\mathcal{I})$,
its  \emph{$\mathcal{S}$-line graph}
(\emph{$\mathcal{S}$-point graph})  is the graph
with vertex set $\mathcal{L}$ (vertex set $\mathcal{P}$) and two
vertices are adjacent if and only if they are incident with a common
point (line, respectively).

\begin{defi}\label{cliq-pls-const-def}
{\rm Let  $\Gamma \in \mathcal{F}(m,r)$ with $m\geq 2,r\geq 1$. Let
$\mathcal{P}=V(\Gamma)$ and $\mathcal{L}=V(C(\Gamma))$, and
$\mathcal{I}\subseteq \mathcal{P}\times \mathcal{L}$ be the set of
pairs $(p,l)$ such that $p\in l$. Let  $\mathcal{S}(\Gamma)$ be the
triple $(\mathcal{P},\mathcal{L},\mathcal{I})$, and let
$\overline{\mathcal{S}(\Gamma)}$ be the graph with vertex set
$\mathcal{P}\cup \mathcal{L}$ and edges all pairs $(p,l)$ such that
$(p,l)\in \mathcal{I}$. }
\end{defi}

By \cite[Lemma 4.7]{DJLP-clique}, we have the following result.

\begin{lemma}\label{cliq-pls-construction}
Let  $\Gamma\in \mathcal{F}(m,r)$ where $m\geq 2,r\geq 1$. Let
$\mathcal{S}(\Gamma)=(\mathcal{P},\mathcal{L},\mathcal{I})$  and $\overline{\mathcal{S}(\Gamma)}$ be as in
Definition \ref{cliq-pls-const-def}. Then  $\mathcal{S}(\Gamma)$ is a partial linear space of order $(m, r + 1)$,  $\overline{\mathcal{S}(\Gamma)}$ is the incidence graph of $\mathcal{S}(\Gamma)$ with girth at least 8. In particular,
$\overline{\mathcal{S}(\Gamma)}$ is a bipartite graph with one bipart corresponding to the vertex set of
$\Gamma$ and the other bipart corresponding to the vertex set of
$C(\Gamma)$.
\end{lemma}

Let  $\Gamma\in \mathcal{F}(m,r)$ where $m\geq 2,r\geq 1$. Let
$\mathcal{S}(\Gamma)=(\mathcal{P},\mathcal{L},\mathcal{I})$  and $\overline{\mathcal{S}(\Gamma)}$ be as in
Definition \ref{cliq-pls-const-def}.    Let
$\textbf{a}=(x_0,x_1,\ldots,x_s)$ be an $s$-arc ($s$-geodesic) of
$\overline{\mathcal{S}(\Gamma)}$. Then $\textbf{a}$ is a
\emph{point-$s$-arc} (\emph{point-$s$-geodesic}) if $x_0$ lies in
$V(\Gamma)$.
Let $G\leq \Aut(\overline{\mathcal{S}(\Gamma)})$. If for any $t\leq
s$ and any two point-$t$-arcs (point-$t$-geodesics)
$\textbf{a},\textbf{b}$ starting from the same vertex $u$ of
$\overline{\mathcal{S}(\Gamma)}$, there exists $g\in G_u$ such that
$\textbf{a}^g=\textbf{b}$, then $\overline{\mathcal{S}(\Gamma)}$ is
said to be \emph{point-$(G,s)$-arc transitive}
(\emph{point-$(G,s)$-geodesic transitive}).

The following result will be used frequently in the remaining of the paper.

\begin{theo}{\rm (\cite{JT-local})}\label{cliq-pls-theo-12}
Let  $\Gamma\in \mathcal{F}(m,r)$ with $m\geq 2,r\geq 1$ and let
$G\leq \Aut(\Gamma)$. Then the following statements are equivalent.

{\rm (1)} $\Gamma$ is  $(G,2)$-geodesic transitive.

{\rm (2)} $\overline{\mathcal{S}(\Gamma)}$ is point-$(G,4)$-arc
transitive and   locally $(G,3)$-arc transitive.

{\rm (3)}   $\overline{\mathcal{S}(\Gamma)}$ is  point-$(G,4)$-arc
transitive and line-$(G,1)$-arc transitive.
\end{theo}

Let $L,R$ be
subgroups of $G$ such that $L\cap R$ is core-free. We define a
\emph{bicoset graph} $\Cos(G,L,R)$ to be a bipartite graph with two
biparts $\Delta_1=[G:L]$ and $\Delta_2=[G:R]$, and edges $\{Lx,Ry\}$
whenever $xy^{-1}\in LR$. Note that $xy^{-1}\in LR$ if and only if
$Lx\cap Ry\neq \emptyset$.

Let $\Sigma \in \mathcal{F}(m,r)$ with $m\geq 2,r\geq 1$ be $(G,2)$-geodesic transitive and let
$\mathcal{S}=\overline{\mathcal{S}(\Sigma)}$ be as in
Definition \ref{cliq-pls-const-def}. Then
$\mathcal{S}$ is biregular with two parts, and  vertices in one part, say $\Delta_1$, have valency $m$
and vertices in the other part, say  $\Delta_2$, have valency $r+1$.

Suppose that $G$ is quasiprimitive of  on $\Delta_1$ but not  quasiprimitive
on $\Delta_2$. Then there exists $1\neq N\unlhd G$ such that $N$ is
not transitive on $\Delta_2$. Further, the quotient graph $\mathcal{S}_N\cong \K_{1,k}$
where $k$ is the number of orbits of $N$ in $\Delta_2$. (Let $\Gamma$ be a graph, and $N\leq \Aut(\Gamma)$. Suppose that
$\mathcal{B}=\{B_1,B_2,\ldots,B_n \}$ is the set of orbits of $N$ on  $V(\Gamma)$.
Then the \emph{quotient graph} $\Gamma_{N}$  is defined to be the graph with vertex set
$\mathcal{B}$ such that $\{B_i,B_j\}$ is an edge of
$\Gamma_{N}$ if and only if there exist $x\in B_i, y\in
B_j$ such that $\{x,y\}$ is an edge of $\Gamma$.)

Moreover, if  $G$ is quasiprimitive of type HS on $\Delta_1$ but not  quasiprimitive
on $\Delta_2$, then
\cite[Theorem 1.1]{GLP-2006} indicates the following result.

\begin{lemma}\label{cliq-noths-01}
Let $\Gamma \in \mathcal{F}(m,r)$ be a
$(G,2)$-geodesic transitive graph where $m\geq 2,r\geq 1$.
Suppose that $G$ is quasiprimitive of type HS on $\Delta_1$ but not  quasiprimitive
on $\Delta_2$, where $\Delta_i$ is as above. Then  $\mathcal{S}\cong \Cos(X,X_v,X_w)$ with
$X=T\times T,X_v=\{(t,t)|t\in T\},X_w=\{(m_1h,m_2h)|m_i\in M,h\in
H\}$ where $M$ is the center of the largest normal $p$-subgroup of
$H$, and $T,H$ are as in one of the lines of  Table \ref{HS},
and $p$ is the defining characteristic of $T$.
\end{lemma}

\begin{table}[t]\caption{HS  }\label{HS}
\medskip
\centering
\begin{tabular}{|c|c|c|c|c|c|}
\hline $T$     &  $H$ & $M$  \\
\hline   $PSL(2,q)$    & $[q]:\mathbb{Z}_{(q-1)/(2,q-1)}$ & $[q]$  \\
\hline  $Sz(2^{l})$    & $[q^2]:\mathbb{Z}_{(q-1)}$, $q=2^l,l\geq 3$ odd & $[q]$  \\
\hline  $Ree(3^{l})'$    & $[q^3]:\mathbb{Z}_{(q-1)}$, $q=3^l,l\geq 3$ odd & $[q]$  \\
\hline $PSU(3,q)$, $q\geq 4$    & $[q^3]:\mathbb{Z}_{(q^2-1)/(3,q+1)}$  & $[q]$  \\
\hline $PSL(d,q)$, $d\geq 3$    & $[q^{d-1}]:(\mathbb{Z}_{q-1}.SL(d-1,q))\mathbb{Z}_{(q-1,d-1)}$ & $[q^{d-1}]$   \\
  \hline
\end{tabular}
\end{table}

\section{ Proof of Theorem \ref{2gt-cliq-basic-2}}

Let $\Gamma \in \mathcal{F}(m,r)$ be a
$(G,2)$-geodesic transitive graph where $m\geq 2,r\geq 1$.
The following lemma shows that  $G$ does not act quasiprimitively
of type $\HS$ on the vertex set.
We use  $[q]$ to denote  a group of order $q$.

\begin{lemma}\label{cliq-noths-1}
Let $\Gamma \in \mathcal{F}(m,r)$ be a
$(G,2)$-geodesic transitive graph where $m\geq 2,r\geq 1$. Then $G$ does not act quasiprimitively
of type $\HS$ on $V(\Gamma)$.
\end{lemma}
\proof Suppose to the contrary that $G$  acts quasiprimitively of type HS on
$V(\Gamma)$.  Since $\Gamma$ is locally isomorphic to $m\K_r$  for some $m\geq
2,r\geq 1$, it follows  from  \cite[Lemma 4.7]{DJLP-clique} that $\Gamma$ is the point graph of a
partial linear space $\mathcal{S}(\Gamma)$ such that the incidence
graph $\Sigma:=\overline{\mathcal{S}(\Gamma)}$ has valency $(m,r+1)$ and
girth at least $ 8$. Further, by Theorem \ref{cliq-pls-theo-12},
\begin{center}
$\Sigma$   is  point-$(G,4)$-arc  transitive  and   locally  $(G,3)$-arc  transitive;
\end{center}
and by \cite[Lemma 4.7]{DJLP-clique},  the $\mathcal{S}(\Gamma)$-line graph  is isomorphic to $C(\Gamma)$. Let
$V(\Sigma)=\Delta_1\cup \Delta_2$ where
$\Delta_1=V(\Gamma)$ and $\Delta_1=V(C(\Gamma))$. As  $\Sigma$
is locally $(G,3)$-arc transitive,   it follows from \cite[Theorems 1.1 and
1.2]{GLP1} that $G$ is not quasiprimitive on $\Delta_2$. Thus there exists
$1\neq N\unlhd G$ which has $k\geq 2$ orbits in $\Delta_2$, and
the quotient graph $\Sigma_N\cong \K_{1,k}$. Therefore, $\Sigma$, $G$ and $N$ satisfy the (STAR) condition of  \cite[p.642]{GLP-2006}.

Let the socle  of $G$ be  $X:=T\times T$ where $T$ is a nonabelian simple group. Then  $T$ is regular on $\Delta_1$, as $G$ acts quasiprimitively of type HS on $\Delta_1$.
Moreover, by  Lemma \ref{cliq-noths-01},    $T$ is one of the groups listed in  Table \ref{HS}.
Since $T$ is a normal subgroup of $ G$ and $C_G(T)=T$, it follows that $G/T\lesssim \Aut(T)$, and so $G\leq T.\Aut(T)$.  Let $L=\{(t,t)|t\in T\}$ and $R=\{(l_1h,l_2h)|l_i\in M,h\in H\}$, where $M$ is the center of the largest normal $p$-subgroup of
$H$. Then $\Sigma\cong \Cos(X,L,R)$.
Let $u=L$, $v=R$ and $w=Lr (\neq L)$. Then $(u,v,w)$ is a 2-geodesic.
Further, $X_{u,v}=L\cap R$ and $X_{u,v,w}\cong \{(x,x)|x\in K\}$
where $K\leq H$ is one of the groups in Table \ref{HS}. In the remainder, we inspect the groups in Table \ref{HS}.
Let $(u,v,w,x,y)$ be a 4-arc of $\Sigma$ and let $A:=T.\Aut(T)$. Then $G\leq A$.

(i)   $T\cong Sz(q)$ where $q=2^l$, $l\geq 3$ odd, $H\cong [q^2]:\mathbb{Z}_{q-1}$ and
$M\cong [q]$. Then $|T|=q^2(q^2+1)(q-1)$ and the valency of $\Sigma$ is $(q^2+1,q)$.
Since
$|A|=|T|^2.|\Out(T)|$,  $|\Out(T)|=l$ and  $T$ is regular on $\Delta_1$, it follows that  $|\Delta_1|=|T|$, and
so  $|A_u|=|T|.|\Out(T)|$.

Since $\Sigma$ is point-$(G,4)$-arc transitive, it is point-$(A,4)$-arc transitive.  The valency of $\Sigma$ is $(q^2+1,q)$, and so
$|A_{u,v}|=q^2(q-1).l$, hence $|A_{u,v,w}|=q^2.l$ and $|A_{u,v,w,x}|=l$.
In particular, $q-1$ divides $l$. However, as  $l\geq 3$ is  odd, $q-1$ does not divide $l$, a contradiction. Thus $T\ncong Sz(q)$ for $q=2^l$.

(ii)   $T\cong PSU(3,q)$ where $q=p^f\geq 4$, $p$ is prime,   $H\cong
[q^3]:\mathbb{Z}_{(q^2-1)/(3,q+1)}$ and $M\cong [q]$.  In this case, $|T|=q^3(q^2-1)(q^3+1)/(3,q+1)$, and the valency of $\Sigma$ is $(q^3+1,q)$. Moreover, by
\cite[p.249]{DM-1}, $H$ has a normal subgroup of order $q^3$ acting
regularly on $\Sigma(u)$.

The order of the group $A$ is  $|A|=|T|^2.|\Out(T)|$ where $|\Out(T)|=2f.(3,q+1)$. Since
$|\Delta_1|=|T|$, it follows that
$$|A_u|=|T|.|\Out(T)|=q^3(q^2-1)(q^3+1).2f,$$
and as
the valency of $\Sigma$ is $(q^3+1,q)$, we have
$$|A_{u,v}|=q^3(q^2-1).2f,$$
hence
$$|A_{u,v,w}|=q^3(q+1).2f$$
and
$$|A_{u,v,w,x}|=(q+1).2f.$$
Since $\Gamma$ is point-$(A,4)$-arc transitive, it follows that $q-1$ divides $(q+1)2f$.
Let $b$ be the largest common divisor of $q-1$ and
$q+1$. As $q+1$ is equal to $q-1+2$, it follows that $b=1$ or 2.

Suppose first that $b=1$. As $q-1$ is a divisor of  $(q+1)2f$, it follows that
$q-1$ divides $2f$. Further, both $q-1$ and $q+1$ are odd integers. Hence
$q=2^f$ and $q-1$ divides $f$. Thus $f=1$ and $q=2$, contradicting
$q\geq 4$.
Suppose next that $b=2$. Then $((q-1)/2,(q+1)/2)=1$. As $q-1$ divides
$(q+1)2f$, it follows that $(q-1)/2$ divides $2f$. As $q-1$ is even,
we know that $q=p^f$ is odd. Thus either $p=3$ and $f=1,2$ or $p=5$ and $f=1$ .
Since $q\geq 4$, it follows that either $p=3$ and $f=2$ or $p=5$ and
$f=1$.
Since $X_{u,v,w}=[q^3]$, the order $|A_{u,v,w}|$ divides $q^3.2f.(3,q+1)$. Hence $|A_{u,v,w,x}|$ divides $2f.(3,q+1)$.
Thus $|\Sigma_4(u)\cap \Sigma(x)|=q-1$ is a divisor of  $(2f,q+1)$. However, for both cases $(p=3,f=2)$ and $(p=5, f=1)$, $q-1$ does not divide  $2f.(3,q+1)$,
which is a contradiction. Thus $T\ncong
PSU(3,q)$ for $q=p^f\geq 4$, $p$ is a prime.

(iii)  $T\cong PSL(d,q)$ where $d\geq 3$, $q=p^f$, $p$ is a prime,  $H\cong
[q^{d-1}]:(\mathbb{Z}_{(q-1)}.SL(d-1,q)).\mathbb{Z}_{(d-1,q-1)}$ and $M\cong
[q^{d-1}]$. In this case, we have
$$|T|=(q^d-1)(q^d-q)\cdots
(q^d-q^{d-1})/(q-1)(d,q-1)$$
$$=q^{d-1}(q^d-1)(q^d-q)\cdots
(q^d-q^{d-2})/(d,q-1)$$
and $\Sigma$ has valency $((q^d-1)/(q-1),q^{d-1})$.

Since $X_{u,v,w}\cong [q^{d-1}]$, we know that $|A_{u,v,w}|$ divides $q^{d-1}.2f.(d,q-1)$. Recall that   the girth of $\Sigma$ is
at least  $ 8$. It follows that
$$|\Sigma_3(u)\cap \Sigma(w)|=(q^d-1/q-1)-1=q^d-q/q-1.$$
Thus $|A_{u,v,w,x}|$ divides
$$q^{d-1}.2f.(d,q-1)/(q^d-q/q-1)=q^{d-1}(q-1).2f.(d,q-1)/(q^d-q).$$
Note that
$A_{u,v,w,x}$ is transitive on $\Sigma_4(u)\cap \Sigma(x)$. As
$|\Sigma_4(u)\cap \Sigma(x)|=q^{d-1}$, it follows that  $q^{d-1}-1$
divides $q^{d-1}(q-1).2f.(d,q-1)/(q^d-q)$. Then we have
$$q^{d-1}-1 \leq q^{d-1}(q-1).2f.(d,q-1)/(q^d-q),$$
that is,
$$(q^d-q)(q^{d-1}-1)\leq q^{d-1}(q-1).2f.(d,q-1).$$
Hence
$$(q^{d-1}-1)^2\leq
q^{d-1}(q-1).2f.(d,q-1)/q\leq q^{d-1}.2f.(d,q-1).$$
Since $q\geq 2$ and $d\geq 3$, it follows that  $q^{d-1}-1\geq q^{d-2}$, and so  $q^{d-1}-1\leq q.2f.(d,q-1)$.
The inequation  $q^{d-1}-q\leq q^{d-1}-1$ indicates  that $q^{d-1}-q\leq q.2f.(d,q-1)$, hence
$q^{d-2}-1\leq 2f.(d,q-1)$. Since   $q\geq 2$ and $d\geq 3$, we know that  $q^{d-2}-1\leq 2f.(d,q-1)$ does not occur. It leads to that
$A_{u,v,w,x}$ is not transitive on
$\Sigma_4(u)\cap \Sigma(x)$, which means  that  $T\ncong PSL(d,q)$.

(iv)  $T\cong
Ree(q)$ where  $q=3^l$, $l\geq 3$ odd, $H\cong [q^3]:\mathbb{Z}_{q-1}$ and
$M\cong [q]$. Then the order of $T$ is $|T|=q^3(q^3+1)(q-1)$ and the valency of $\Sigma$ is $(q^3+1,q)$. Note that
$|A|=|T|^2.|\Out(T)|$ and $|\Out(T)|=l$. Since $|\Delta_1|=|T|$, it
follows that $|A_u|=|T|.|\Out(T)|$.
Since $\Sigma$ is point-$(A,4)$-arc transitive and the valency of $\Sigma$ is $(q^3+1,q)$, it follows that
$|A_{u,v}|=q^3(q-1).l$, $|A_{u,v,w}|=q^3.l$ and $|A_{u,v,w,x}|=l$.
In particular, $q-1$ divides $l$.  However, the fact  that $l\geq 3$ is an odd integer indicates that $q-1$ does not divide $l$, a contradiction. Thus
$T$ is not isomorphic to $ Ree(q)$.

(v)   $T\cong
PSL(2,q)$ where  $q=p^f\geq 4$, $p$ is prime,  $H\cong
[q]:\mathbb{Z}_{(q-1)/(2,q-1)}$ and $M\cong [q]$. In this case  $|T|=q(q+1)(q-1)/(2,q-1)$ and the valency of $\Sigma$ is $(q+1,q)$.
Moreover, the order of $A$ is  $|A|=|T|^2.|\Out(T)|$ and
$|\Out(T)|=2f.(2,q-1)$. Since $|\Delta_1|=|T|$, it follows that
$|A_u|=|T|.|\Out(T)|$.

Again, since $\Sigma$ is a point-$(A,4)$-arc transitive graph with  valency  $(q+1,q)$, it follows that
$$|A_{u,v}|=q(q-1).2f.(2,q-1)/(2,q-1)=q(q-1).2f,$$
\begin{center}
$|A_{u,v,w}|=q.2f$
and $|A_{u,v,w,x}|=2f$.
\end{center}

In particular,  $q-1$ divides $2f$.
If $q-1$ is an odd integer, then $q-1$ divides $2f$, and hence $q-1$ divides
$f$, which is impossible.
If $q-1$ is an even integer, then $q$
is odd and $q-1$ divides $2f$. Hence $p$ is an odd integer, which is again  a contradiction. Hence
$T\ncong PSL(2,q)$.

Therefore, $G$ does not act quasiprimitively
of type HS on $V(\Gamma)$. \qed

Now we  prove the  main theorem.

\medskip

\noindent {\bf Proof of Theorem \ref{2gt-cliq-basic-2}.}
Let $\Gamma \in \mathcal{F}(m,r)$  be a $(G,2)$-geodesic transitive graph  where $m,r\geq 2$ and $G\leq \Aut(\Gamma)$.   Let $\mathcal{S}(\Gamma)$ be the partial linear space
as in Definition \ref{cliq-pls-const-def} and
$\overline{\mathcal{S}(\Gamma)}$ be its incidence graph. Then  it follows from
Theorem \ref{cliq-pls-theo-12}   that
$\overline{\mathcal{S}(\Gamma)}$ is locally $(G,3)$-arc transitive.

Suppose that $G$ acts quasiprimitively on
$V(\Gamma)$   of type $X$. Then  $G$ is quasiprimitive on
the point part of $\overline{\mathcal{S}(\Gamma)}$, and  by Theorems
1.2 and 1.3 of \cite{GLP1}, $X$ is one of the types in the set $\{\HA,\AS,\TW,\PA,\SD,\HS\}$.
Moreover,
Lemma \ref{cliq-noths-1} says that  $X$ is not of the HS type, and so $X\in \{\HA,\AS,$ $\TW,\PA,\SD\}$.
If   the induced action of $G$ on $V(C(\Gamma))$ is not
quasiprimitive, then by \cite[Theorem 1.3]{GLP1},  $X$ is one of
$\{\HA,\PA,\AS,\TW\}$, see Table \ref{Star}; and if the induced action of $G$ on $V(C(\Gamma))$ is  quasiprimitive of type
$Y$, then by \cite[Theorem 1.2]{GLP1}, $(X,Y)$ is one of
$\{(\HA,\HA),(\PA,\PA),(\AS,\AS),$ $(\TW,\TW),$ $(\SD,\PA),(\PA,\SD)\}$, see Table \ref{Biprimitive}. We conclude the proof. \qed

\section{ Proof of Theorem \ref{2gt-basic-th2}}

We recall the definition of Hamming graphs. The \emph{Hamming graph}  $\Gamma=\H(d,n)$  has
vertex set $\Delta^d=\{(x_1,\cdots,x_d)|x_i\in \Delta\}$, the
cartesian product of $d$-copies of $\Delta$, where
$\Delta=\{1,a,\ldots,a^{n-1}\}$, $d\geq 2$  and $n\geq 2$.  Then two
vertices are adjacent if and only if they are different
in exactly one coordinate.

\begin{examp}\label{2gt-qp-examphamm1}
{\rm  Let $\Gamma=\H(d,n)$ where $d\geq 2$ and $n\geq 5$. Then  $\Gamma \in \mathcal{F}(d,n-1)$,  and by \cite[Proposition 2.2]{DJLP-compare}, $\Gamma$ is $(G,2)$-geodesic transitive for $G=\Aut(\Gamma)\cong S_n\wr S_d$.   In particular,  $G$ acts primitively of
type PA on $V(\Gamma)$. Let $N=S_n^d\lhd G$. Then $N$ has $d$ orbits on $V(C(\Gamma))$, and each orbit being the set of
maximal cliques of $\Gamma$ associated with a given coordinate. Hence $G$ is not quasiprimitive on  $V(C(\Gamma))$.
 }
\end{examp}

\medskip
Let $\Gamma=(V(\Gamma),E(\Gamma))$ be a $(G,2)$-geodesic transitive graph. Suppose that $G$ is quasiprimitive on $V(\Gamma)$ of Product Action type.  Then $G$ preserves a Cartesian decomposition $V(\Gamma)=\Delta^k=\Delta_1\times \cdots \times \Delta_k$. Set $\Delta=\{\delta_1=\delta,\delta_2,\ldots,\delta_m\}$. Since $G$ is transitive on $V(\Gamma)$, it follows that $G$ induces a natural transitive action on $\Omega=\{1,2,\ldots,k\}$. Assume that
the  induced transitive group on $\Omega$ by $G$ is $P$, that is, $P\leq \Sym(\Omega)$. On the other hand, the set stabilizer $G_{\Delta}$ induces a faithful quasiprimitive permutation  group $H$ on $\Delta$ of type Almost Simple, that is, $H\leq \Sym(\Delta)$. Thus, $G\leq H\wr P$.

Let $u=(\delta,\delta,\ldots,\delta)$ and let $\Gamma(u)=\{v\in V(\Gamma)|\{u,v\}\in E(\Gamma)\}$.
Let $\Gamma\circ \Gamma(u)=\{v\in V(\Gamma)|\Gamma(u)\cap \Gamma(v)\neq \emptyset\}$. Then $\Gamma_2(u)\subseteq \Gamma\circ \Gamma(u)\subseteq \Gamma(u)\cup \Gamma_2(u)$.

\begin{lemma}\label{pa-distincts}
Let $\Gamma$ be a $(G,2)$-geodesic transitive graph. Suppose that $G$ is primitive on $V(\Gamma)$ of PA type and preserves a Cartesian decomposition $V(\Gamma)=\Delta^k=\Delta_1\times \cdots \times \Delta_k$. Let $X=soc(H)$ be transitive on $\Delta$.  Let $u=(\delta,\delta,\ldots,\delta)$ and  $v=(v_1,\ldots,v_k)\in \Gamma(u)$. Then  $u,v$ have $j$ distinct coordinates where  $j=1$ or $2$.
\end{lemma}
\proof  Suppose that $u,v$ have $j$ distinct coordinates where $j\geq 2$.  Since $G$ is transitive on the set $\{\Delta_1,\Delta_2,\ldots,\Delta_k\}$ and $u=(\delta,\delta,\ldots,\delta)$, it follows that $G_u$ induces a transitive action on $\{1,2,\ldots,k\}$. We  assume that $v_1\neq \delta$ and $v_2\neq \delta$.

Suppose that $(\delta,v_1)\in \theta_1$ and $(\delta,v_2)\in \theta_2$ where $\theta_1,\theta_2$ are two orbits of $X$ in $\Delta\times \Delta$. Then
$(v_1,\delta)\in \theta_1^*$ and $(v_2,\delta)\in \theta_2^*$. Since $G$ is primitive on $V(\Gamma)$ of Product Action type,
$X=soc(H)$ is primitive  on $\Delta$ of Almost Simple type.

Let $\Omega=\{w |w\in \Delta,w^x=w,\forall x\in X_{v_1}\}$.
Since $H$ is primitive on $\Delta$ and  $X$ is a normal subgroup of $H$,
it follows that $X$ is transitive on $\Delta$.
Suppose that $\Omega\cap \Omega^g\neq \emptyset$ for some $g\in H$, and say
$r\in \Omega\cap \Omega^g$.
Then $r^g=t\in \Omega$, hence $X_r^g=X_{r^g}=X_t$. Since  $r,t\in \Omega$, we have
$$X_r=X_{v_1r}=X_{v_1}=X_{v_1t}=X_t,$$
and so
$$X_r^g=X_r.$$
Note that
$\{w |w\in \Delta,w^x=w,\forall x\in X_{v_1}^g\}=\Omega^g$. Thus
$\Omega=\Omega^g$, and so $\Omega$ is a block of $H$. Note that  $H$ is primitive on $\Delta$, we must have $|\Omega|=1$, that is, $\Omega=\{v_1\}$.

Similarly,
$\{r|r^x=r,\forall x\in X_{v_2}\}=\{v_2\}$. Thus,
there exist $x\in X_{v_1}$ and $y\in X_{v_2}$ such that $\delta^x\neq \delta$,  and $\delta^y\neq \delta$.
Denoting $\delta^x=\alpha$ and  $\delta^y=\beta$.
Then  $(v_1,\alpha)=(v_1^x,\delta^x)=(v_1,\delta)^x\in \theta_1^*$ and $(v_2,\beta)=(v_2^y,\delta^y)=(v_2,\delta)^y\in \theta_2^*$.
In particular,  $\overline{x}=(x,1,1,\ldots,1)\in G$, $\overline{y}=(1,y,1,\ldots,1)\in G$ and $\overline{z}=(x,y,1,\ldots,1)\in G$.

Let $\overline{\alpha}=(\alpha,\delta, \dots,\delta)$, $\overline{\beta}=(\delta,\beta,\delta,\dots,\delta)$ and $\overline{\gamma}=(\alpha,\beta,\delta,\dots,\delta)$ be three vertices.
Since $v=(v_1,\ldots,v_k)\in \Gamma(u)$, it follows that
$$(v,u)^{\overline{x}}=(v,\overline{\alpha}), (v,u)^{\overline{y}}=(v,\overline{\beta}), (v,u)^{\overline{z}}=(v,\overline{\gamma}).$$
Thus $v\in \Gamma(u)\cap \Gamma(\overline{\alpha})\cap \Gamma(\overline{\beta})\cap \Gamma(\overline{\gamma})$.

Since $\Gamma$ is $(G,2)$-geodesic transitive, $G_u$ is transitive on $\Gamma(u)$.
Since $v\in \Gamma(u)$ and $u,v$ have $j$ distinct coordinates where $j\geq 2$, it follows that $\overline{\alpha},\overline{\beta}\in \Gamma\circ \Gamma(u)\setminus \Gamma(u)=\Gamma_2(u)$. Again, since  $\Gamma$ is $(G,2)$-geodesic transitive, $G_u$ is transitive on $\Gamma_2(u)$, it follows that $\overline{\gamma} \notin \Gamma_2(u)$, so
$\overline{\gamma} \in \Gamma(u)$. Thus, $j=2$.
\qed

Let $\Delta=\{0,1,2,\ldots,m-1\}$ and $\Delta^k=\Delta\times \cdots \times \Delta$ where $m,k\geq 2$. Define  $\Gamma$ to be the graph with vertex set $\Delta^k$, and two vertices $u=(u_1,\ldots,u_k)$ and $v=(v_1,\ldots,v_k)$ are adjacent if and only if they have  exactly $2$ different coordinates.

\begin{lemma}\label{2distinct}
Let $\Gamma$ be a graph defined as above. If   $\Gamma$ is $2$-geodesic  transitive, then $k=2$ and $m\geq 3$ and  $\Gamma$ is isomorphic to $\overline{\H(2,m)}$.
\end{lemma}
\proof Suppose that $\Gamma$ is  $2$-geodesic  transitive. If $k=2$ and $m=2$, then vertices $(0,0)$ and $(1,0)$ are not in the same connected component, and so $\Gamma$ is disconnected.

If $k=2$ and $m\geq 3$, then $\Gamma$ is connected with diameter 2. Moreover, $\Gamma$ is isomorphic to $\overline{\H(2,m)}$.

Assume that $k= 3$. If $m=2$,
then the vertex $(0,0,0)$ lies in a connected component with  $4$ vertices, and $(1,0,0)$ lies in another connected component with also $4$ vertices, and so  $\Gamma$ is disconnected.
Suppose that  $m\geq 3$. Let  $u=(0,0,0)$, $v_1=(1,1,0)$, $w_1=(2,0,0)$, $v_2=(2,1,0)$, and $w_2=(1,1,1)$. Then $(u,v_1,w_1)$
and $(u,v_2,w_2)$ are two $2$-geodesics. Noting that the stabilizer of $u$ in the automorphism group  can not map $w_1$ to $w_2$, contradicts $\Gamma$ is $2$-geodesic  transitive.

Assume that $k\geq 4$. Suppose that  $m\geq 3$.
 Let  $u=(0,0,0,0,0^{k-4})$, $v_1=(1,1,0,0,0^{k-4})$, $w_1=(1,2,0,1,0^{k-4})$,  and $w_2=(1,1,1,1,0^{k-4})$. Then $(u,v_1,w_1)$
and $(u,v_1,w_2)$ are two $2$-geodesics. Noting that the stabilizer of $u$ in the automorphism group  can not map $w_1$ to $w_2$, contradicts that $\Gamma$ is $2$-geodesic  transitive.  Hence $m=2$.
However, in this case, vertex $(0,0,\ldots,0)$ lies in a connected component with  $2^{k-1}$ vertices, and $(1,0,\ldots,0)$ lies in another connected component with also $2^{k-1}$ vertices, and so  $\Gamma$ is disconnected.
  \qed

\begin{lemma}\label{pa-distincts-loc}
Let $\Gamma$ be a $(G,2)$-geodesic transitive graph. Suppose that $G$ is primitive on $V(\Gamma)$ of PA type and preserves a Cartesian decomposition $V(\Gamma)=\Delta^k=\Delta_1\times \cdots \times \Delta_k$. Let $X=soc(H)$ be transitive on $\Delta$. Suppose that $fix_{\Delta}(X_{\delta})=\{\delta\}$. Let $u=(\delta,\delta,\ldots,\delta)$ and  $v=(v_1,\ldots,v_k)\in \Gamma(u)$.  Assume that  $v_1$ is in the orbit $\Sigma_i(\delta)$ of $H_\delta$ in $\Delta$ and the orbital graph $(\Delta,\Sigma_i)$ is connected. If $u,v$ have exactly $1$ distinct coordinate, then $H$ is $2$-transitive on $\Delta$.
\end{lemma}
\proof  Suppose that   $u,v$ have exactly $1$ distinct coordinate.
Say  $v=(v_1\neq \delta,\delta,\ldots,\delta)$.   Assume that  $v_1$ is in the orbit $\Sigma_i(\delta)$ of $H_\delta$ in $\Delta$.
Since $H$ is a  primitive permutation  group  on $\Delta$,
 the orbital graph $(\Delta,\Sigma_i)$ is connected.
Let $\alpha \in \Sigma_i(\delta)$. Then there exists $x\in H_\delta$ such that $v_1^x=\alpha$. Thus $\overline{x}=(x,1,\ldots,1)\in G$
and $(u,v)^{\overline{x}}=(u,(\alpha,\delta,\ldots,\delta))$.

Let $\beta \in \Sigma_{i^*}(\alpha)$. Then  $(\alpha,\beta)\in \Sigma_{i^*}$, and so
$(\beta,\alpha)\in \Sigma_i$.  Recall that $\alpha \in \Sigma_i(\delta)$, and so $\alpha \in \Sigma_i(\beta)\cap \Sigma_i(\delta)$,
hence in the orbital graph $(\Delta,\Sigma_i)$, we have $\beta \in \Sigma_i\circ \Sigma_i(\delta)$.

Let $w=(\alpha,\beta,\delta,\ldots,\delta)$. Since $(\delta,\alpha)\in \Sigma_i$ and $(\beta,\alpha)\in \Sigma_i$, there exists $y\in H_\alpha$
such that $\delta^y=\beta$. Then
$$(1,y,1,\ldots,1)\in G$$
and
$$(u,(\alpha,\delta,\ldots,\delta))^{(1,y,1,\ldots,1)}=((\delta,\beta,\delta,\ldots,\delta),w).$$
On the other hand, $(y,1,1,\ldots,1)\in G$ and $$(u,(\alpha,\delta,\ldots,\delta))^{(y,1,1,\ldots,1)}=((\beta,\delta,\ldots,\delta),(\alpha,\delta,\ldots,\delta)).$$
Hence $(\beta,\delta,\ldots,\delta)$ is in $\Gamma(u)\cup \Gamma_2(u)$, as $(\alpha,\delta,\ldots,\delta)\in \Gamma(u)$. Assume that $(\beta,\delta,\ldots,\delta)$ is in $ \Gamma_2(u)$.
Then since $\Gamma$ is $(G,2)$-geodesic transitive, $G_u$ has an element that maps  $(\beta,\delta,\ldots,\delta)$ to $w$, which is impossible.
Thus $(\beta,\delta,\ldots,\delta)$ is in $ \Gamma(u)$.
Note that $G_u$ can map $(\delta,\beta,\delta,\ldots,\delta)$ to $(\beta,\delta,\ldots,\delta)$.
Thus $(\delta,\beta,\delta,\ldots,\delta)$ is in $ \Gamma(u)$.
Recall that $((\delta,\beta,\delta,\ldots,\delta),w)$ is an arc, it follows that $w$ is in $\Gamma_2(u)$.
Further, $\{(\alpha,\beta,\delta,\ldots,\delta)|\alpha \in \Sigma_i(\delta),\beta \in \Sigma_{i^*}(\alpha)\}$ is contained in $\Gamma\circ \Gamma(u)\setminus \Gamma(u) $.

Let  $r\in \Sigma_i\circ \Sigma_i(\delta)$. Then  $\Sigma_i(\delta)\cap \Sigma_i(r)\neq \emptyset$, and say $r'\in \Sigma_i(\delta)\cap \Sigma_i(r)$.
Then $(\delta,r')\in \Sigma_i$ and $(r,r')\in \Sigma_i$, and so there exists $x''\in H_{r'}$ such that $\delta^{x''}=r$.
Recall that   $v_1\in \Sigma_i(\delta)$.
Then there exists $x^{'''}\in H_\delta$ such that $v_1^{x^{'''}}=r'$. Thus $\overline{x}=(x^{'''},1,\ldots,1)\in G$
and $(u,v)^{\overline{x}}=(u,(r',\delta,\ldots,\delta))$. Further, $\tilde{x}=(x'',1,\ldots,1)\in G$ and $(u,(r',\delta,\ldots,\delta))^{\tilde{x}}=(\overline{r},(r',\delta,\ldots,\delta))$,
where $\overline{r}=(r,\delta,\ldots,\delta)$. Thus
we have  $\overline{r}\in \Gamma\circ \Gamma(u)$. Since  $G_u$ is transitive on $\Gamma_2(u)$ and $w\in \Gamma\circ \Gamma(u)\setminus \Gamma(u)=\Gamma_2(u)$, it follows that $\overline{r}\in \Gamma(u)$.
Thus, there exists $\overline{g}\in G_u$ such that $(u,v)^{\overline{g}}=(u,\overline{r})$.  Further, $\overline{g}$ induces an element $g\in H_\delta$
such that $v_1^{g}=r$. Hence  $r\in \Sigma_i(\delta)$, and so
$\Sigma_i\circ \Sigma_i(\delta)\subseteq \Sigma_i(\delta)$.

Since  the orbital graph $(\Delta,\Sigma_i)$ is connected, for any
$z\in \Delta\setminus  \{\delta\}$, there is a path $(\delta=w_1,w_2,w_3\ldots,w_n=z)$, where $(w_i,w_{i+1})\in \Sigma_i$. Hence  $w_2\in \Sigma_i(w_1)\cap \Sigma_i(w_3)$,  it leads to   $w_3\in \Sigma_i\circ\Sigma_i(w_1)$. By the previous argument,
$\Sigma_i\circ\Sigma_i(w_1) \subseteq \Sigma_i(w_1)$, we have that  $w_3\in \Sigma_i(w_1)$. Thus   $w_3\in \Sigma_i(w_1)\cap \Sigma_i(w_4)$,  so  $w_4\in \Sigma_i\circ\Sigma_i(w_1) \subseteq \Sigma_i(w_1)$. By induction, for each $i=2,\ldots,n$,  $w_i\in \Sigma_i\circ\Sigma_i(w_1) \subseteq \Sigma_i(w_1)$. Thus
$\Delta \setminus \{w_1\}\subseteq  \Sigma_i(w_1)$, and so   $\Delta\subseteq \{w_1\}\cup \Sigma_i(w_1)$, hence $H$ is 2-transitive on $\Delta$.
\qed

For every vertex $u$ of $\Gamma$, we define  $\Gamma\circ \Gamma(u)=\{v\in V(\Gamma)|\Gamma(u)\cap \Gamma(v)\neq \emptyset\}$.
Then $\Gamma\circ \Gamma(u)\setminus \Gamma(u)=\Gamma_2(u)$.

We are ready to   prove the second theorem.

\medskip
\noindent {\bf Proof of Theorem \ref{2gt-basic-th2}.}
Let $\Gamma$ be a $2$-geodesic transitive graph. Suppose that
$A:=\Aut(\Gamma)$ acts primitively on $V(\Gamma)$ of PA type. Then
$A$ preserves a Cartesian decomposition $V(\Gamma)=\Delta^k$.
Let $H$ be the induced subgroup of $A_{\Delta}$ in $\Delta$. Then since  $A$ is primitive on $V(\Gamma)$ of PA type,    $H$ is primitive on $\Delta$. Let  $u\in V(\Gamma)$. The graph  $\Gamma$ is $2$-geodesic transitive implying that  both $\Gamma(u)$ and $\Gamma\circ \Gamma(u)\setminus \Gamma(u)=\Gamma_2(u)$ are orbits of $A_u$ in $V(\Gamma)\setminus \{u\}$.
Let  $v\in \Gamma(u)$. Then it follows from  Lemma \ref{pa-distincts} that  $u$ and $v$ have $j$ distinct coordinates where $j=1,2$.
If $j=2$, then   by   Lemma \ref{2distinct},
$k=2$ and $n=|\Delta|\geq 3$ and  $\Gamma$ is isomorphic to $\overline{\H(2,n)}$.
If $n=3$ or 4, then $G$ is primitive on the vertex set of type HA,  a contradiction. Thus $n\geq 5$.

Assume that   $j=1$. Then    Lemma \ref{pa-distincts-loc} says that $H$ is 2-transitive on $\Delta$. Hence  $v\in \Gamma(u)$ if and only if  $u,v$ have exactly $1$ distinct coordinate.
In this case  $\Gamma$ is the  Hamming graph $\H(d,n)$ where $d\geq 2$ and $n\geq 2$.
If $n= 2$, then  $\Gamma$ is both antipodal and bipartite, a contradiction. If $n=3$ or 4, then $G$ is primitive on the vertex set of type HA, again a contradiction. Thus $n\geq 5$.
  \qed

\begin{lemma}\label{2gt-star-pa}
Let $\Gamma \in \mathcal{F}(m,r)$ for some  $m,r\geq 2$ be a $(G,2)$-geodesic transitive graph. Suppose that $G$ is
quasiprimitive on $V(\Gamma)$ of type \PA.  If $G$ is
not quasiprimitive on $V(C(\Gamma))$, then $\Gamma$ is a Hamming graph and $G$ is in \cite[Example 5.1]{GLP-2006}.
\end{lemma}
\proof Let $\overline{\mathcal{S}(\Gamma)}$ be as in Definition \ref{cliq-pls-const-def}. Since $\Gamma \in \mathcal{F}(m,r)$ is
$(G,2)$-geodesic transitive, it follows from Theorem \ref{cliq-pls-theo-12} that   $\overline{\mathcal{S}(\Gamma)}$  is
locally $(G,3)$-arc transitive.  Since $G$ is
quasiprimitive on $V(\Gamma)$ but not quasiprimitive on $V(C(\Gamma))$, $\overline{\mathcal{S}(\Gamma)}$ satisfies the condition (STAR) of \cite[p.642]{GLP-2006}. Then by
\cite[Theorem 1.2]{GLP-2006}, $\Gamma$ is a Hamming graph and $G$ is in \cite[Example 5.1]{GLP-2006}. \qed

\section{Examples}

In this section, we give some examples for the quasiprimitive action types AS, HA and SD, but for quasiprimitive action
type TW, we still do not know if it can occur or not. We are investigating this problem at the moment.

\subsection{Almost Simple}

For a group $G$ and a subgroup $H$, denote by $[G:H]$ the set of right cosets of $H$ in
$G$. Let $D=HgH \subset G$ with $g\notin H$ and $g^2\in H$. Then the \emph{coset graph}
$\Gamma=\Cos(G,H,D)$ has vertex set $[G:H]$ and arc set $\{(Hg,Hdg)|g\in G,d\in D\}$, and
$\Gamma$ is connected if and only  if $G=\langle D\rangle$,
$\Gamma$ is undirected if and only  if $D=D^{-1}$,  its valency is $|D:H|$.

\begin{examp}\label{smallval6-psl-1}
{\rm  Let $G=PSL(2,p)$ where $p$ is a prime,  $p\equiv \pm 1\pmod{24}$.
Choose  $h\in H$ with $o(h)=4$ where $H\cong S_4$ is a maximal
subgroup of $G$.  Let $g\in C_G(\langle h\rangle)\setminus H$ with
$o(g)=3$. Let $\Gamma=\Cos(G,H,HgH)$. }
\end{examp}

\begin{lemma}\label{smallval6-psl-1-1}
Let groups $G$, $H$,  elements $h$, $g$ and graph  $\Gamma=\Cos(G,H,HgH)$ be as in Example \ref{smallval6-psl-1}. Then

{\rm (1)} $N_G(\langle h\rangle) =C_G(\langle h\rangle):\langle
z\rangle$ where $o(z)=2$ and $z\in H$,  $G=\langle H,g\rangle$ and
$HgH=Hg^{-1}H$;

{\rm (2)} $\Gamma$ has girth  $3$ and $|H:H^g\cap H|=6$;

{\rm (3)}   $[\Gamma(u)]\cong 3\K_2$,  $\Gamma$ is $(G,2)$-geodesic
transitive, and  $G$ acts quasiprimitively of $\AS$ type on both $V(\Gamma)$ and $V(C(\Gamma))$.
\end{lemma}
\proof  (1) Since $H$ is a maximal subgroup of $G$ and
$g\notin H$, it follows  that $G=\langle H,g\rangle$.

Since $h$ is an element of order 4 in $H$, there exists an
involution  $z\in H$ such that $h^z=h^{-1}$. Thus $z\in N_G(\langle
h\rangle)\setminus C_G(\langle h\rangle)$. Since $C_G(\langle
h\rangle)\cong \mathbb{Z}_{ \frac{p\mp1}{2}} $ and $N_G(\langle
h\rangle)\cong D_{ p\mp 1}$, it follows that $N_G(\langle h\rangle)
=C_G(\langle h\rangle):\langle z\rangle$. Let $C_G(\langle
h\rangle)=\langle x\rangle $. Then $x^z=x^{-1}$, and hence
$g^z=g^{-1}$. Thus $HgH=HzgzH=Hg^{-1}H$.

(2) Since $(H,Hg,Hg^2)$ is a triangle, it follows that $\Gamma$ has girth
$3$.

Since $z\in H$, it follows that $HgH=HgzH$. Let $f=gz$. Then
$\Gamma=\Cos(G,H,HgH)$ $=\Cos(G,H,HfH)$, $o(f)=2$ and $G=\langle
H,f\rangle$. Further the valency of $\Gamma$ is $|H:H^f\cap H|$. Let
$K=\langle h\rangle$. Since $h^f=h^{gz}=h^{-1}$, it follows that
$K^f=K$, and hence $K\leq H^f\cap H$. Now we prove that $K= H^f\cap
H$.

Since $H^f\cap H$ is a subgroup of $H$ and $K\leq H^f\cap H$, it
follows that $H^f\cap H=K,J\cong D_8$ or $H$. If $H^f\cap H=H$, then
$H$ is a normal subgroup of $G$, contradicting the fact that $G$ is
a nonabelian simple group. Suppose that  $H^f\cap H=J$.  Since
$z^2=1$ and $h^z=h^{-1}$, it follows that $\langle h,z\rangle\cong
D_8$. Since the subgroup of $H$ which contains $h$ and isomorphic to
$D_8$ is unique, it follows that $J=\langle h,z\rangle$. Thus
$z^f\in H$. However, $z^f=z^{gz}=zg\notin H$, a contradiction. Thus
$H^f\cap H=K$, and hence $|H:H^f\cap H|=6$.

(3) Let $(u,v)$ be an arc where $u=H$ and $v=Hf$. Then
$(u,v)^f=(v,u)$.  Let $E$ be the kernel of the $G_u$-action on
$\Gamma(u)$. Then $E\trianglelefteq G_u$. Since $\Gamma$ is $G$-arc
transitive, it follows that $6\big | |G_u^{\Gamma(u)}|$, and hence
$|E|\leq 4$. Suppose that $E\neq 1$. Then $E\cong \mathbb{Z}_2\times \mathbb{Z}_2$.
Further, $E\leq G_{u,v}$. Since $|G_{u,v}|=4$, it follows that $E=
G_{u,v}$, and hence $G_{u,v}\trianglelefteq G_u$. Further, since
$(u,v)^f=(v,u)$, it follows that $G_{u,v}^f=G_{u,v}$. Since
$G=\langle G_u,f\rangle$, $G_{u,v}\trianglelefteq G$ contradicting
that $G$ is nonabelian simple. Thus, $E=1$.

Since $\Gamma$ is not complete, $|\Gamma(u)\cap \Gamma(v)|=1,2,3,4$.
If $|\Gamma(u)\cap \Gamma(v)|=4$, then $\Gamma(u) \cong \K_{3[2]}$,
and hence $\Gamma \cong \K_{4[2]}$, a contradiction. Suppose that
$|\Gamma(u)\cap \Gamma(v)|=3$. If $[\Gamma(u)]$ does not contain
triangles, then $\Gamma\cong \K_{3[3]}$, a contradiction. Assume that
$[\Gamma(u)]$ contains a triangle. Suppose that $\Gamma(u)=\{v_i|i
=1,2,3,4,5,6\}$ where $v_1=v$. Then $(v_1,v_2,v_3)$ and
$(v_4,v_5,v_6)$ are triangles, and $v_j$ is adjacent with $v_{j+3}$
where $j\in \{1,2,3\}$, that is, $[\Gamma(u)]$ is the triangular
prism. Thus $G_u^{\Gamma(u)} \lesssim S_3\times \mathbb{Z}_2$, a
contradiction. If $|\Gamma(u)\cap \Gamma(v)|=2$, then
$[\Gamma(u)]\cong C_6$ or $2C_3$.  Since $G_u^{\Gamma(u)}\cong
G_u\cong S_4$ is not isomorphic to a subgroup of $D_{12}$, it
follows that  $[\Gamma(u)]\ncong C_6$. If $[\Gamma(u)]\cong 2C_3$,
then $G_u\cong S_4$ has two imprimitive blocks on $\Gamma(u)$ which
is also impossible. Thus, $|\Gamma(u)\cap \Gamma(v)|=1$, that is,
$[\Gamma(u)]\cong 3\K_2$.

Since $E=1$ and  $|G_{u,v}|=4$, it follows that every element of
$G_{u,v}$ moves some vertices of $\Gamma_2(u)\cap \Gamma(v)$. It
follows that $G_{u,v}$ is regular on $\Gamma_2(u)\cap \Gamma(v)$.
Thus $\Gamma$ is $(G,2)$-geodesic transitive. Finally, since $G$ is a nonabelian simple group,  $G$ acts quasiprimitively of AS type on both $V(\Gamma)$ and $V(C(\Gamma))$. \qed

\begin{examp}\label{2gt-qp-examp1}
{\rm (1)  Let $\Omega=\{1,2,3,4,5,6\}$. Let  $\Gamma$ be the graph with vertex set the set of 2-subsets of $\Omega$, and
two vertices are adjacent if and only if they are disjoint. Then  $\Gamma\cong KG(6,2) \in \mathcal{F}(3,2)$, of diameter $2$ and
$\Aut(\Gamma)\cong S_6$.  It is known that  $\Gamma$ is  arc transitive.
 Let $(u,v)$ be an arc where $u=\{1,2\},v=\{3,4\}$. Then $\Gamma_2(u)\cap \Gamma(v)=\{\{1,5\},\{1,6\},\{2,5\},\{2,6\}\}$.
Hence $A_{u,v}$ is transitive on  $\Gamma_2(u)\cap \Gamma(v)$, and so $\Gamma$ is $2$-geodesic transitive.
Further, $\Aut(\Gamma)$ acts primitively of type AS on both $V(\Gamma)$ and $V(C(\Gamma))$.

(2)  Let $\Gamma$ be the Foster graph. Then by  \cite[p.221]{BCN} and \cite[Theorem
1.1]{Weiss-2},  $\Gamma$ is bipartite $5$-arc transitive with two parts $\Delta_1$ and $\Delta_2$ and has 90 vertices, both girth and
diameter are $8$, and valency is 3.

Let $\Sigma$ be the graph with vertex set $\Delta_1$ or $\Delta_2$, and two vertices are  adjacent in $\Sigma$ if and only if they are
at distance 2 in $\Gamma$.  Then
$|V(\Sigma)|=45$, and $\Sigma\in \mathcal{F}(3,2)$ is 3-fold
cover of $KG(6,2)$.  Since $\Gamma$ is 5-arc transitive, it follows from Theorem \ref{cliq-pls-theo-12}
that $\Sigma$ is $2$-geodesic transitive. Further,  the induced action of $\Aut(\Gamma)$ on $V(\Sigma)$ is
quasiprimitive of type AS.
 }
\end{examp}

\begin{examp}\label{primitive-lemma-1}
{\rm

(1) Let $G,H,k,s$ be as in \cite[Table 1]{Li-2001-biprimitive}, and let $\Sigma=\Cos(G,H,HgH)$ where $g\in G\setminus H$ with $g^2\in H$. Then by
\cite[Theorem 1.3]{Li-2001-biprimitive}, $\Sigma$ is $(G,s)$-transitive of valency $k$ with $s\geq 4$ and $G=\Aut(\Sigma)$.

Let $\Gamma$ be the standard double cover of $\Sigma$. Then by \cite[Theorem 1.4 (i)]{Li-2001-biprimitive}, $G$ acts  primitively of type AS on each bipart vertex set of $\Gamma$, and $\Gamma$ is $t$-transitive of valency $k$ with $t\geq 4$. Thus,  either $\diam(\Gamma)\leq 3 $ or $\Gamma$ has girth at least $8$.
We prove that $\Gamma$ has girth at least $8$. Suppose that $\diam(\Gamma)\leq 3 $.
Let $u$ and $u'$ be an antipodal vertex pair. Then the distance between them is at most $3$.
Suppose that   the distance between $u$ and $u'$ is $r$. Let  $(u,\ldots,u')$ be an $r$-arc from $u$ to $u'$ and $(u,\ldots,z)$ be an $r$-arc  where $z\neq u'$. Since $u,u'$ are in the same antipodal block, it follows that $G$ can not map $(u,\ldots,u')$ to $(u,\ldots,z)$,  contradicts that $\Gamma$ is $t$-transitive with $t\geq 4$. Hence $\Gamma$ has girth at least $8$, and so,  $\Gamma$ is a partial linear space. Thus, its point graph and line graph are both in $\mathcal{F}(k,k-1)$ and vertex primitive of type AS. Further, by Theorem \ref{cliq-pls-theo-12}, they are 2-geodesic transitive.

(2)  Let $G,H,k,s$ be as in Table \ref{biquasi-as}, and let $\Gamma=\Cos(G,H,HgH)$ where $g\in G\setminus H$ with $g^2\in H$.  Then by \cite[Theorem 1.4 (ii)]{Li-2001-biprimitive}, $\Gamma$ is a bipartite  $s$-transitive graph of valency $k$ and $\Aut(\Gamma)$ acts primitively of type AS on each vertex bipart. Since  $\Gamma$ is a bipartite  $s$-transitive with $s\geq 5$, and since $s\leq (g+2)/2$ where $g$ is the girth of $\Gamma$, it follows that the girth of $\Gamma$ is at least $8$, and so  $\Gamma$ is the incidence graph of a partial linear space. Therefore, the  point graph and the line graph of this partial linear space are both in $\mathcal{F}(k,k-1)$,
and vertex primitive of type AS. Further, by Theorem \ref{cliq-pls-theo-12}, they are 2-geodesic transitive.

\begin{table}[]\caption{ }\label{biquasi-as}
\medskip
\centering
\begin{tabular}{|c|c|c|c|c|c|}
\hline $T=soc(G)$     &  $G$ & $H\cap T$ & $k$ & $s$  \\
\hline   $Sp(4,2^m)$    & $G\leq \Aut(T)$ &  $[2^{3m}]:GL(2,2^m)$ & $2^m+1$  & $5$    \\
\hline   $G_2(3^m)$    & $G\leq \Aut(T)$ &  $[3^{5m}]:GL(2,3^m)$ & $3^m+1$  & $7$    \\
  \hline
\end{tabular}
\end{table}

}

\end{examp}

We use $E(\Gamma)$ to denote the edge set of $\Gamma$.

\begin{examp}\label{subdivision-as-star-1}
{\rm   Let $T=PSL(3,q)$ for prime power $q$ and $G=\Aut(T)$. Let
$\Delta_1$ be the set of $1$-dimensional subspaces of a
$3$-dimensional vector space over finite field $GF(q)$ and let $\Delta_2$ be the set of $2$-dimensional
subspaces. We define $\Gamma$ to be the bipartite graph with vertex
set $\Delta_1\cup \Delta_2$ with adjacency given by inclusion. Then
$G\leq \Aut(\Gamma)$. By \cite[Example 2.2]{GL-2010-edge},  $\Gamma$ is $(G,4)$-arc
transitive and $G$ is primitive of type AS on
$E(\Gamma)$. Since $\Gamma$ is bipartite, $G$ is not quasiprimitive on $V(\Gamma)$.
Let the valency of $\Gamma$ be $m$. Let $\Sigma=C(\Gamma)=L(\Gamma)$. Then by \cite[Theorem 1.4]{DJLP-clique}, $\Sigma\in
\mathcal{F}(2,m-1)$, and by \cite[Theorem 1.3]{DJLP-line}, $\Sigma$ is $(G,2)$-geodesic transitive. Further, $G$ is primitive of type AS on
$V(\Sigma)$ but not quasiprimitive on the set of maximal cliques of
$\Sigma$ (that is, $V(\Gamma)$).

If  $(d,r,q)=(3,1,2)$, the graph $\Gamma$ is the Heawood
graph.

}
\end{examp}

\begin{examp}\label{2gt-star-as}
{\rm Let $\Gamma \in \mathcal{F}(m,r)$  be a $(G,2)$-geodesic transitive graph where $m,r\geq 2$.
Let $\mathcal{S}(\Gamma)=(\mathcal{P},\mathcal{L},\mathcal{I})$  and $\overline{\mathcal{S}(\Gamma)}$ be as in
Definition \ref{cliq-pls-const-def}. Then $\Gamma$ is isomorphic to the point graph of $\mathcal{S}(\Gamma)$, and  by Theorem \ref{cliq-pls-theo-12} that $\overline{\mathcal{S}(\Gamma)}$ is   locally $(G,3)$-arc transitive.

Suppose that $G$ is
quasiprimitive on $V(\Gamma)$ of type AS with socle $T$.  If $G$ is
not quasiprimitive on $V(C(\Gamma))$, then $T$ has $k\geq 2$ orbits on $V(C(\Gamma))$. In particular,
$\overline{\mathcal{S}(\Gamma)}$, $G$ and $T$ satisfy \cite[p.642, (STAR)]{GLP-2006}.
Thus by \cite[Theorem 1.3]{GLP-2006}, $T$ and $k$ lie in  \cite[Table 1.2]{GLP-2006}.

%Further, if $T\cong E_6(p^f)$ or  $^{2}E_6(p^6)$, we don't have examples of  $\overline{\mathcal{S}(\Gamma)}$, but $\overline{\mathcal{S}(\Gamma)}$
%exists for the other  $T$, see \cite[p.643-644]{GLP-2006}. Thus
%if $T\cong E_6(p^f)$ or  $^{2}E_6(p^6)$, we don't have examples of  $\Gamma$, but $\Gamma$ exists for the other  $T$.
}
\end{examp}

\subsection{Holomorph Affine}

%\begin{lemma}\label{2gt-ha-lem1}
%Let $\Gamma \in \mathcal{F}(m,r)$  be a $(G,2)$-geodesic transitive graph  where $m,r\geq 2$ and $G\leq \Aut(\Gamma)$. Suppose that $G$ is
%quasiprimitive on $V(\Gamma)$ of type HA.  Then
%$\Gamma$ is a Cayley graph, say $\Cay(T,S)$ over the socle $T$ of $G$ and either

% \begin{itemize}
%\item[(i)] $[S]\cong r\K_2$ and $T\cong \mathbb{Z}_3^{r'}$ where $2\leq r'\leq r$; or

%\item[(ii)]  $[S]\cong r\K_{2^n-1}$ and $T\cong \mathbb{Z}_2^{r'n}$ where $2\leq r'\leq r$.
%\end{itemize}
%\end{lemma}

\begin{constr}{\rm (\cite[Construction 4.1]{GLP-2006})}\label{star-ha}
{\rm Let $V$ be a $d$-dimensional vector space over $GF(q)$ where $q=p^f$
for some prime $p$, and let $G_0$ be a subgroup of $GL(d,q)$. Let $\mathcal{N}$ be a collection of $m$-dimensional
subspaces of $V$ permuted $2$-transitively by $G_0$ such that $\mathcal{N}$ spans $V$, and
suppose that for $W\in \mathcal{N}$, $(G_0)_W$ acts transitively on the set of non-trivial elements
of $W$. Let $\mathcal{B}$ be the set of translates of the elements of $\mathcal{N}$, $G=V:G_0$ and $\Sigma$ be the
incidence graph relative to $V$ and $\mathcal{B}$.}
\end{constr}

\begin{examp}\label{2gt-star-ha}
{\rm Let $\Gamma \in \mathcal{F}(m,r)$  be a $(G,2)$-geodesic transitive graph where $m,r\geq 2$.  Let
$\mathcal{S}(\Gamma)=(\mathcal{P},\mathcal{L},\mathcal{I})$  and $\overline{\mathcal{S}(\Gamma)}$ be as in
Definition \ref{cliq-pls-const-def}. Then $\Gamma$ is isomorphic to the point graph of $\mathcal{S}(\Gamma)$, and  by Theorem \ref{cliq-pls-theo-12} that $\overline{\mathcal{S}(\Gamma)}$ is   locally $(G,3)$-arc transitive.

Suppose that $G$ is quasiprimitive on $V(\Gamma)$ of type HA but
not quasiprimitive on $V(C(\Gamma))$. Then $G$ is quasiprimitive on the point graph of $\mathcal{S}(\Gamma)$ of type HA but
not quasiprimitive on the line graph of $\mathcal{S}(\Gamma)$. Thus $(\overline{\mathcal{S}(\Gamma)},G)$ satisfies  \cite[p.642, (STAR)]{GLP-2006}, and  by \cite[Theorem 4.2]{GLP-2006}, $\overline{\mathcal{S}(\Gamma)}$ is the graph in  Construction \ref{star-ha}.
}
\end{examp}

\begin{examp}\label{2gt-qp-examphamm12}
{\rm   Let $\Gamma=\H(d,n)$ where $d\geq 2$ and $n\in
\{3,4\}$. Then $\Gamma \in \mathcal{F}(d,n-1)$, and by \cite[Proposition 2.2]{DJLP-compare}, $\Gamma$ is $(G,2)$-geodesic transitive
for $G=\Aut(\Gamma)$.  In particular,  $G$ acts primitively of type HA on
$V(\Gamma)$. Let $N=S_n^d\lhd G$. Then $N$ has $d$ orbits on $V(C(\Gamma))$, and each orbit being the set of maximal
cliques of $\Gamma$ associated with a given coordinate. Hence $G$ is not quasiprimitive on  $V(C(\Gamma))$.
 }
\end{examp}

\begin{examp}\label{2gt-localconnect-priexampha}
{\rm (1) Let $n\geq 7$ be odd and $V=F_2^n$ be an $n$-dimensional
vector space over field $F_2$. Let $W=\{(w_1,\ldots,w_n)\in
V|w_1+\cdots+w_n=0\}$. Let $\Gamma$ be a graph with vertex set $W$,
and two vectors are adjacent if they differ in precisely 2 entries.
Then $\Gamma$ is the half-cube, and it is 2-geodesic transitive with
diameter $[\frac{n}{2}]$, valency $\frac{n(n-1)}{2}$,
$|V(\Gamma)|=2^{n-1}$, see \cite[p.264]{BCN}. Further, $\Gamma$ is
locally connected and $\Aut(\Gamma)\cong S_2^{n-1}. S_n$ acts
primitively of type HA on $V(\Gamma)$.

{\rm (2)} Let $n\geq 5$ be odd and $V=F_2^n$ be an $n$-dimensional
vector space over field $F_2$. Let $W=\{(w_1,\ldots,w_n)\in
V|w_1+\cdots+w_n=0\}$. Let $\Gamma$ be a graph with vertex set $W$,
and two vectors are adjacent if they agree in  precisely 1 entry.
Then $\Gamma$ is the folded  cube, and it is 2-geodesic transitive
with diameter $[\frac{n}{2}]$, valency $\frac{n(n-1)}{2}$,
$|V(\Gamma)|=2^{n-1}$, see \cite[p.264]{BCN}. Further,
$girth(\Gamma)\geq 4$ and $\Aut(\Gamma)\cong S_2^{n-1}. S_n$ acts
primitively of type HA on $V(\Gamma)$.

{\rm (3)} Let $n\geq 3$ be even and $V=F_2^n$ be an $n$-dimensional
vector space over field $F_2$. Let $W=\{(w_1,\ldots,w_n)\in
V|w_1+\cdots+w_n=0\}$. Let $X\subseteq W$ and $U=W/X$.  Let $\Gamma$
be a graph with vertex set $U$, and two vectors $w_1+X$, $w_2+X$ are
adjacent iff $w_1,w_2$ either differ in  precisely 2 entries, or
agree in precisely 2 entries. Then $\Gamma$ is the folded half-cube,
and it is 2-geodesic transitive with $|V(\Gamma)|=2^{n-2}$, see
\cite[p.264]{BCN}. Further, $girth(\Gamma)\geq 4$ and
$\Aut(\Gamma)\cong S_2^{n-2}. S_n$ acts primitively of type HA on
$V(\Gamma)$.
}
\end{examp}

\subsection{Simple Diagonal}

\begin{examp}\label{2gt-cliq-sd-pa}
{\rm Let $\Gamma=\Cos(G,L,R)$ be the bicoset graph in \cite[Theorem1.1
(1)]{GLP-2007-5arc}. Then $\Gamma$ is locally $5$-arc transitive but
not locally $6$-arc transitive of valency $(2^m+1,2^m)$ where $m$ is a positive integer. Let $A=\Aut(\Gamma)$, let $\Delta_1$ be
the part of valency $2^m+1$ and $\Delta_2$ be the part of valency
$2^m$. Then $A$ acts primitively of type SD on $\Delta_1$ and
primitively of type PA on $\Delta_2$.

It's known that a locally $s$-arc transitive graph of girth $g$ with at least one vertex of valency greater than $2$ satisfies $s\leq (g+2)/2$. Since $\Gamma$ is locally $5$-arc transitive, it follows that the girth $g$ of $\Gamma$ satisfies
$g\geq 8$, and so $\Gamma$ is the incidence graph of a partial linear space $\mathcal{S}$. Let $\Sigma$ be the  point graph of $\mathcal{S}$  and  $\Sigma'$ be the  line graph  of $\mathcal{S}$. Then by \cite[Theorem 1.4]{DJLP-clique}, $\Sigma=C(\Sigma')$ and $\Sigma'=C(\Sigma)$, and either  $\Sigma$ or $\Sigma'$ lies in  $\mathcal{F}(2^m+1,2^m-1)$,  the other one lies in $\mathcal{F}(2^m,2^m)$.
Since $\Gamma$ is locally $5$-arc transitive, it follows from Theorem \ref{cliq-pls-theo-12} that both $\Sigma$ and $\Sigma'$  are 2-geodesic transitive.

If $\Sigma$  lies in $\mathcal{F}(2^m+1,2^m-1)$,  then $\Sigma'$ lies in $\mathcal{F}(2^m,2^m)$, and  $V(\Sigma)=\Delta_1$, $V(\Sigma')=\Delta_2$. Thus by \cite[Theorem1.1(3)]{GLP-2007-5arc},  $A$ acts primitively of type SD on $V(\Sigma)$ and
primitively of type PA on $V(\Sigma')=V(C(\Sigma))$. If $\Sigma'$  lies in $\mathcal{F}(2^m+1,2^m-1)$,  then $\Sigma$ lies in $\mathcal{F}(2^m,2^m)$, and  $V(\Sigma')=\Delta_1$, $V(\Sigma)=\Delta_2$. Thus by \cite[Theorem1.1(3)]{GLP-2007-5arc}, $A$ acts primitively of type SD on $V(\Sigma')$ and
primitively of type PA on $V(\Sigma)=V(C(\Sigma'))$.

}
\end{examp}

\end{document}